\documentclass[10pt]{article}
\usepackage{amsmath, amssymb, braket, amsthm, graphicx, titlefoot, color}
\usepackage{setspace}
\usepackage[export]{adjustbox}
\usepackage{float}
\usepackage[labelfont=bf]{caption}
\usepackage[square,numbers]{natbib}

\title{Path Planning in a Riemannian Manifold using Optimal Control}
\author{Souma Mazumdar\thanks{Email: \texttt{souma.mazumdar@bose.res.in},  Phone: 09903144810} \\ Department of Theoretical Sciences \\ S. N. Bose National Centre for Basic Sciences \\ Block - JD, Sector - III, Salt Lake City, Kolkata - 700 106 }
\date{}

\begin{document}
	\maketitle 
	
	\begin{abstract}
       We consider the motion planning of an object in a Riemannian manifold where the object is steered from an initial point to a final point utilizing optimal control. Considering Pontryagin Minimization Principle we compute the Optimal Controls needed for steering the object from initial to final point. The Optimal Controls were solved with respect to time $t$ and shown to have norm $1$ which should be the case when the extremal trajectories, which are the solutions of Pontryagin Principle, are arc length parametrized. The extremal trajectories are supposed to be the geodesics on the Riemannian manifold. So we compute the geodesic curvature and the Gaussian curvature of the Riemannian structure. 
       
	\end{abstract}

	\section{Introduction}
    Path Planning or Motion Planning is today an active area of research in the control theory community which concerns about steering an object from one point to another utilizing controls. Practical application of this is abound in real life scenario. Starting from Optimal path planning of a mobile robot in an unknown terrain to proper manouevering of an autonomous vehicle to landing of an aircraft in an airport. Utilizing the Optimal Control Theory and Pontryagin Principle is one of the tools used by Control theorists to deal with motion planning problems. While on the other hand there are other techniques employed in motion planning problems like the continuation method as shown by Sussman. All of these methods concerns about finding a proper control which can steer the object from one point to another. Sussman in his work\cite{sussmann1992new} discussed about a number of techniques of motion planning all of which are different in the mathematical ideas they involve. One of them is the Continuation method where the mathematical idea goes like lifting the problem of motion planning from the state space to the control space and then solve a differential equation to find the controls. A number of interesting papers\cite{chitour2006continuation,chitour1998line,sussmann1993continuation} have been written along this line by Sussman and Chitour. We will not go far discussing about various methods of motion planning but come back to the topic of this paper which is motion planning in a Riemannian manifold utilizing optimal controls.\par
    Most of the problems in control theory fall under the purview of Sub-Riemannian Geometry. While speaking about Sub-Riemannian Geometry we think of a underlying distribution $\mathcal{D}$ which is a subspace of the the tangent space $T_{q}M$ at a point $q$, where $q \in M$, $M$ being the manifold. Here the vector fields span a subspace of the tangent space but if the vector fields form a bracket generating family then the vector fields together with their Lie brackets span the whole tangent space. But in our problem we speak of Riemannian manifold. The prime difference in the Riemannian case from the Sub-Riemannian problem is that the vector fields span the whole tangent space. We deal with the problem of motion planning in a Riemannian manifold using the principles of Optimal Control. Moreover the manifold being Riemannian, it has an intrinsic curvature. So Optimal Control problems in a Riemannian manifold will be very different from other cases. In Sub-Riemannian case the extremal curves which are the solutions of Pontryagin Principle can be both Normal and Abnormal extremals. But in Riemannian case as shown by Agrachev\cite{agrachev2019comprehensive} the extremals are always Normal extremals. So in our work we are primarily concerned with Normal extremals as Abnormal extremals do not exist in a Riemannian manifiold.\par  
    In our problem we consider the manifold to be a $2-D$ Riemannian manifold. So the control problem has two orthogonal vector fields which span the tangent space $T_qM$ at a point $q \in M$. The family of vector fields can be involutive, that is their Lie bracket generates $0$ or may be non involutive. We apply the principles of Optimal Control theory and Pontryagin Principle to find the Optimal Controls which help in steering the object from one point to another on the manifold. Further we assume, the normal extremals which are the solutions of Pontryagin principle are arc length parametrized. So the velocity at any point on the trajectory is always $1$. As the manifold possesses an intrinsic curvature, in addition to computing the Optimal Controls, we find the geodesic and Gaussian curvatures. The Gaussian curvature evaluates to $0$ which should be the case in case of involutive vector fields\cite{agrachev2019comprehensive} as the Lie bracket of our chosen family of orthogonal vector fields generates $0$.
    \section{Problem Formulation}
    We consider a $2-D$ Riemannian manifold. The vectors fields span the whole tangent space. So in our problem we consider two vector fields which form an orthogonal frame. A control system on such a manifold is described by the following equation.
    \begin{equation}
    \begin{split}
    \dot{q}=u_{1}\begin{bmatrix}
    q_{1} \\ q_{2}
    \end{bmatrix}+u_{2}\begin{bmatrix}
    q_{2} \\ -q_{1}
    \end{bmatrix}
    \end{split}
    \end{equation}
    where $u_{1}$ and $u_{2}$ are the controls. Note that the vector fields $\begin{bmatrix}
    q_{1} \\ q_{2}
    \end{bmatrix}$ and $\begin{bmatrix}
    	q_{2} \\ -q_{1}
    \end{bmatrix}$
    form an orthogonal frame. We cast the problem in the Pontryagin Principle formalism and solve for the Optimal Controls. Pontryagin Principle is a Hamiltonian formalism. So it generates differential equations of position$(q)$ and momentum$(p)$. In our case also we are confronted with such differential equations which are nonlinear and hard to solve. But with the help of the fact that the velocity at any point of the extremal curves is always $1$ we were able to solve these equations with ingenuity.
    \section{Pontryagin Principle and Hamiltonian Formalism}
     We consider spacial optimality and time to vary between $0$ and $1$ that is $t \in (0,1)$. The objective function to be optimized is the length functional $\int_{0}^{1}\sqrt{u_{1}^2+u_{2}^2}$ which is same as optimizing the energy functional $\frac{1}{2}\int_{0}^{1}(u_{1}^2+u_{2}^2)$. We cast the problem is the Pontryagin Formalism. The associated Hamiltonian is constructed as follows
     \begin{equation}
     \begin{split}
     H(p,q,u)=-\frac{1}{2}(u_{1}^2+u_{2}^2)+p^{T}(u_{1}f_{1}(q)+f_{2}(q))
     \end{split}
     \end{equation}
     where $f_{1}(q)=\begin{bmatrix}
     q_{1} \\ q_{2}
     \end{bmatrix}$ and $f_{2}(q)=\begin{bmatrix}
     q_{2} \\ -q_{1}
     \end{bmatrix}$ and $p^{T}$ stands for the transpose of the momentum covector. Here $q=\begin{bmatrix}
     q_{1} \\ q_{2}
        \end{bmatrix}$ and $p=\begin{bmatrix} p_{1} \\ p_{2} \end{bmatrix}$. \\
     Minimizing the above Hamiltonian with respect to the controls we obtain the controls as
     \begin{equation}
     \begin{split}\label{Optimal Contols}
     & u_{1}=<p,f_{1}(q)> \\ &
     u_{2}=<p,f_{2}(q)>
     \end{split}
     \end{equation}
     where $<.,.>$ defines the inner product between the vectors and the covectors. The Hamilton's equations are given by the following
     \begin{equation}\label{state eq}
     \begin{split}
     \dot{q}=<p,f_{1}(q)>f_{1}(q)+<p,f_{2}(q)>f_{2}(q)
     \end{split}
     \end{equation} 
     \begin{equation}\label{costate eq}
     \begin{split}
     \dot{p}=-<p,f_{1}(q)><p,D_{q}f_{1}(q)>-<p,f_{2}(q)><p,D_{q}f_{2}(q)>
     \end{split}
     \end{equation}
     where $D_{q}$ stands for the Jacobian. \\
     The momentum covector written as above is a row vector. When written in column vector form it takes the form,
     \begin{equation}
     \begin{split}
     \dot{p}=-<p,f_{1}(q)><D_{q}f_{1}(q),p>-<p,f_{2}(q)><D_{q}f_{2}(q),p>
     \end{split}
     \end{equation}
     When written out in components representation (\ref{state eq}) takes the form,
     \begin{equation}
     \begin{split}
     & \dot{q}_{1}=(p_{1}q_{1}+p_{2}q_{2})q_{1}+(p_{1}q_{2}-p_{2}q_{1})q_{2} \\ &
     \dot{q}_{2}=(p_{1}q_{1}+p_{2}q_{2})q_{2}-(p_{1}q_{2}-p_{2}q_{1})q_{1}
     \end{split}
     \end{equation}
     For finding the equation for momentum covectors we compute the respective Jacobians.
     Now,
     \begin{equation}
     \begin{split}
     D_{q}f_{1}(q)p=\begin{bmatrix}
     1 & 0 \\ 0 & 1
     \end{bmatrix}\begin{bmatrix}
     p_{1} \\ p_{2}
     \end{bmatrix}=\begin{bmatrix}
     p_{1} \\ p_{2}
     \end{bmatrix}
     \end{split}
     \end{equation}
     and
      \begin{equation}
     \begin{split}
     D_{q}f_{2}(q)p=\begin{bmatrix}
     0 & 1 \\ -1 & 0
     \end{bmatrix}\begin{bmatrix}
     p_{1} \\ p_{2}
     \end{bmatrix}=\begin{bmatrix}
     p_{2} \\ -p_{1}
     \end{bmatrix}
     \end{split}
     \end{equation}
     Therefore (\ref{costate eq}) when written out in components representation takes the form,
     \begin{equation}
     \begin{split}
     & \dot{p}_{1}=-(p_{1}q_{1}+p_{2}q_{2})p_{1}-(p_{1}q_{2}-p_{2}q_{1})p_{2} \\ &
     \dot{p}_{2}=-(p_{1}q_{1}+p_{2}q_{2})p_{2}+(p_{1}q_{2}-p_{2}q_{1})p_{1}
     \end{split}
     \end{equation}
     When simplifield the $4$ differential equations of Hamilton's equation can be written as follows
     \begin{equation}\label{Hamilton eqns}
     \begin{split}
     & \dot{q}_{1}=p_{1}(q_{1}^2+q_{2}^2) \\ &
     \dot{q}_{2}=p_{2}(q_{1}^2+q_{2}^2) \\ &
     \dot{p}_{1}=-q_{1}(p_{1}^2-p_{2}^2)-2p_{1}p_{2}q_{2} \\ &
     \dot{p}_{2}=-q_{2}(p_{2}^2-p_{1}^2)-2p_{1}p_{2}q_{1}
     \end{split}
     \end{equation}
     Moreover we consider the normal extremals are arc length parametrized. Therefore the velocity at any point of the trajectory is always $1$. So we have the additional equation,
     \begin{equation}
     \begin{split}
     \dot{q}_{1}^2+\dot{q}_{2}^2=1
     \end{split}
     \end{equation}
     \section{Solve of the Hamilton's Equations}
     From the $1st$ and $2nd$ equations of the set(\ref{Hamilton eqns}) we have
     \begin{equation}
     \begin{split}
     & p_{1}=\frac{\dot{q}_{1}}{q_{1}^2+q_{2}^2} \\ &
     p_{2}=\frac{\dot{q}_{2}}{q_{1}^2+q_{2}^2}
     \end{split}
     \end{equation}
     Differentaiting the above equations with respect to time we have,
     \begin{equation}
     \begin{split}
     & \dot{p}_{1}=\frac{(q_{1}^2+q_{2}^2)\ddot{q}_{1}-\dot{q}_{1}(2q_{1}\dot{q}_{1}+2q_{2}\dot{q}_{2})}{(q_{1}^2+q_{2}^2)^2} \\ &
     \dot{p}_{2}=\frac{(q_{1}^2+q_{2}^2)\ddot{q}_{2}-\dot{q}_{2}(2q_{1}\dot{q}_{1}+2q_{2}\dot{q}_{2})}{(q_{1}^2+q_{2}^2)^2}
     \end{split}
     \end{equation} 
     Puting the values of $\dot{p}_{1}$ and $\dot{p}_{2}$ in the $3rd$ and $4th$ equations of the set (\ref{Hamilton eqns}) and simplifying we have,
     \begin{equation}\label{doubledot q1}
     \begin{split}
     (q_{1}^2+q_{2}^2)\ddot{q}_{1}-\dot{q}_{1}(2q_{1}\dot{q}_{1}+2q_{2}\dot{q}_{2})=-q_{1}(\dot{q}_{1}^2-\dot{q}_{2}^2)-2\dot{q}_{1}\dot{q}_{2}q_{2}
     \end{split}
     \end{equation}
     \begin{equation}\label{doubledot q2}
     \begin{split}
     (q_{1}^2+q_{2}^2)\ddot{q}_{2}-\dot{q}_{2}(2\dot{q}_{1}q_{1}+2\dot{q}_{2}q_{2})=-q_{2}(\dot{q}_{2}^2-\dot{q}_{1}^2)-2\dot{q}_{1}\dot{q}_{2}q_{1}
     \end{split}
     \end{equation}
     Multiplying equation (\ref{doubledot q1}) by $\dot{q}_{1}$ and equation (\ref{doubledot q2}) by $\dot{q}_{2}$ and adding them we have,
     \begin{equation}\label{add ddot q1 and ddot d2}
     \begin{split}
     & (q_{1}^2+q_{2}^2)(\dot{q}_{1}\ddot{q}_{1}+\dot{q}_{2}\ddot{q}_{2})-(2q_{1}\dot{q}_{1}+2q_{2}\dot{q}_{2})(\dot{q}_{1}^2+\dot{q}_{2}^2)\\ & =-\dot{q}_{1}q_{1}(\dot{q}_{1}^2-\dot{q}_{2}^2)-\dot{q}_{2}q_{2}(\dot{q}_{2}^2-\dot{q}_{1}^2)-2\dot{q}_{1}^2\dot{q}_{2}q_{2}-2\dot{q}_{2}^2\dot{q}_{1}q_{1}
     \end{split}
     \end{equation}
     Now we have the additional equation,
     \begin{equation}\label{velocity eq}
     \dot{q}_{1}^2+\dot{q}^2=1
     \end{equation}
     Differentiating the above equation with respect to time and dividing by $2$ we have
     \begin{equation}\label{derivative velocity eq}
     \dot{q}_{1}\ddot{q}_{1}+\dot{q}_{2}\ddot{q}_{2}=0
     \end{equation}
     Now using equations (\ref{velocity eq}) and (\ref{derivative velocity eq}) and puting in equation (\ref{add ddot q1 and ddot d2}) we have,
     \begin{equation}
     \begin{split}
     -(2q_{1}\dot{q}_{1}+2q_{2}\dot{q}_{2})=-\dot{q}_{1}q_{1}(\dot{q}_{1}^2-\dot{q}_{2}^2)-\dot{q}_{2}q_{2}(\dot{q}_{2}^2-\dot{q}_{1}^2)-2\dot{q}_{1}^2\dot{q}_{2}q_{2}-2\dot{q}_{2}^2\dot{q}_{1}q_{1}
     \end{split}
     \end{equation}
     Simplifying the above equation we have,
     \begin{equation}
     \begin{split}
     -2q_{1}\dot{q}_{1}-2q_{2}\dot{q}_{2}=-(\dot{q}_{1}q_{1}+\dot{q}_{2}q_{2})(\dot{q}_{1}^2+\dot{q}_{2}^2)
      \end{split}
     \end{equation}
     Again using (\ref{velocity eq}) we have the above equation as
     \begin{equation}
     -2q_{1}\dot{q}_{1}-2q_{2}\dot{q}_{2}=-\dot{q}_{1}q_{1}-\dot{q}_{2}q_{2}
     \end{equation}
     which simplifies to 
     \begin{equation}\label{derivative postion eq}
     \dot{q}_{1}q_{1}+\dot{q}_{2}q_{2}=0
     \end{equation}
     Integrating the above equation with respect to time we have,
     \begin{equation}
     q_{1}^2+q_{2}^2=K
     \end{equation}
     where $K$ is the constant of integration. The above constant can be scaled to $1$. Then we have the above equation as,
     \begin{equation}\label{position eqn}
     q_{1}^2+q_{2}^2=1
     \end{equation}
     Using this relation for the $1st$ and $2nd$ equations of the set (\ref{Hamilton eqns}) we have,
     \begin{equation}\label{momentum wrt position}
     \begin{split}
     & \dot{q}_{1}=p_{1} \\ &
     \dot{q}_{2}=p_{2}
     \end{split}
     \end{equation}
     Differentiating the above set with respect to time we have,
     \begin{equation}\label{derivative of momentum wrt to position}
     \begin{split}
     &\dot{p}_{1}=\ddot{q}_{1} \\ &
     \dot{p}_{2}=\ddot{q}_{2}
     \end{split}
     \end{equation}
     Now using equations (\ref{momentum wrt position}) and (\ref{derivative of momentum wrt to position}) and puting in the $3rd$ and $4th$ equations of the set (\ref{Hamilton eqns}) we have,
     \begin{equation}\label{double dot q1 eqn wrt q1}
     \ddot{q}_{1}=-q_{1}(\dot{q}_{1}^2-\dot{q}_{2}^2)-2\dot{q}_{1}\dot{q}_{2}q_{2}
     \end{equation}
     \begin{equation}\label{double dot q2 eqn wrt q2}
     \ddot{q}_{2}=-q_{2}(\dot{q}_{2}^2-\dot{q}_{1}^2)-2\dot{q}_{1}\dot{q}_{2}q_{1}
     \end{equation}
     Now multiplying equation (\ref{double dot q1 eqn wrt q1}) by $q_{2}$ and equation (\ref{double dot q2 eqn wrt q2}) by $q_{1}$ and adding them and then simplifying we have,
     \begin{equation}
     q_{2}\ddot{q}_{1}+q_{1}\ddot{q}_{2}=-2\dot{q}_{1}\dot{q}_{2}(q_{1}^2+q_{2}^2)
     \end{equation}
     Now using equation (\ref{position eqn}) we have the above equation as
     \begin{equation}
     q_{2}\ddot{q}_{1}+q_{1}\ddot{q}_{2}=-2\dot{q}_{1}\dot{q}_{2}
     \end{equation}
     Dividing both sides by $\dot{q}_{1}\dot{q}_{2}$ we have the above  equation as,
     \begin{equation}\label{right side devoid of q}
     \frac{q_{2}\ddot{q}_{1}}{\dot{q}_{1}\dot{q}_{2}}+\frac{q_{1}\ddot{q}_{2}}{\dot{q}_{1}\dot{q}_{2}}=-2
     \end{equation}
     Now from equation (\ref{derivative velocity eq}) we have 
     \begin{equation}
     \frac{\ddot{q}_{1}}{\dot{q}_{2}}=-\frac{\ddot{q}_{2}}{\dot{q}_{1}}
     \end{equation}
     Using the above equation we have equation (\ref{right side devoid of q}) as
     \begin{equation}
     \frac{q_{2}}{\dot{q}_{1}}\left(-\frac{\ddot{q}_{2}}{\dot{q}_{1}}\right)+\frac{q_{1}\ddot{q}_{2}}{\dot{q}_{1}\dot{q}_{2}}=-2
     \end{equation}
     which on simplifying gives,
     \begin{equation}\label{right side devoid of q 2nd eq}
     \frac{\ddot{q}_{2}}{\dot{q}_{1}}\left(\frac{q_{1}}{\dot{q}_{2}}-\frac{q_{2}}{\dot{q}_{1}}\right)=-2
     \end{equation}
     Now from equation (\ref{derivative postion eq}) we have,
     \begin{equation}
     \frac{q_{1}}{\dot{q}_{2}}=-\frac{q_{2}}{\dot{q}_{1}}
      \end{equation}
      Using the above equation and puting in equation (\ref{right side devoid of q 2nd eq}) we have
      \begin{equation}
      2\frac{\ddot{q}_{2}}{\dot{q}_{1}}\frac{q_{1}}{\dot{q}_{2}}=-2
      \end{equation}
      which on simplifying gives
      \begin{equation}
      \frac{\ddot{q}_{2}}{\dot{q}_{2}}=-\frac{\dot{q}_{1}}{q_{1}}
      \end{equation}
      Now the above equation can be written as,
      \begin{equation}
      \frac{d}{dt}(\ln{\dot{q}_{2}})=-\frac{d}{dt}(\ln{q_{1}})
      \end{equation}
      which on integrating gives,
      \begin{equation}\label{derivative q2 wrt q1}
      \begin{split}
      & \ln{\dot{q}_{2}}=-\ln{q_{1}}+\ln{C_{1}} \\ &
      \implies \dot{q}_{2}=\frac{C_{1}}{q_{1}}
      \end{split}
      \end{equation}
      where $C_{1}$ is the constant of integration. \\
      Puting the above value of $\dot{q}_{2}$ in the equation (\ref{velocity eq}) we have,
      \begin{equation}
      \dot{q}_{1}^2+\frac{C_{1}^2}{q_{1}^2}=1
      \end{equation}
      Therefore,
      \begin{equation}
      \begin{split}
      & \dot{q}_{1}=\sqrt{1-\frac{C_{1}^2}{q_{1}^2}} \\ &
      \implies \frac{dq_{1}}{dt}=\sqrt{1-\frac{C_{1}^2}{q_{1}^2}}
      \end{split}
      \end{equation}
      Therefore,
      \begin{equation}
      \begin{split}
      & \frac{dq_{1}}{\sqrt{1-\frac{C_{1}^2}{q_{1}^2}}}=dt \\ &
      \implies \frac{q_{1} dq_{1}}{\sqrt{q_{1}^2-C_{1}^2}} = dt
      \end{split}
      \end{equation}
      Integrating the above equation we have,
      \begin{equation}
      \begin{split}
      & \sqrt{q_{1}^2-C_{1}^2}=t+C_{2} \\ &
      \implies q_{1}=\sqrt{(t+C_{2})^2+C_{1}^2}
      \end{split}
      \end{equation}
      where $C_{2}$ is the constant of integration. \\
      From equation (\ref{derivative q2 wrt q1}) we have,
      \begin{equation}\label{to intergrate q2}
      \begin{split}
      & \dot{q}_{2}=\frac{C_{1}}{\sqrt{(t+C_{2})^2+C_{1}^2}} \\ &
     \implies dq_{2}=\frac{C_{1}dt}{\sqrt{(t+C_{2})^2+C_{1}^2}}
      \end{split}
      \end{equation}
     Now for integration of the right side of the above equation we substitute $t+C_{2}=C_{1}\tan\theta$.
      Therefore by this substitution (\ref{to intergrate q2}) takes the form,
      \begin{equation}
      \begin{split}
      dq_{2}=C_{1}\sec\theta d\theta
      \end{split}
      \end{equation} 
     Integrating both sides we have,
     \begin{equation}
     \begin{split}
     q_{2}=C_{1}[\ln(\sec\theta+\tan\theta)]+C_{3}
     \end{split}
     \end{equation}
     where $C_{3}$ is the constant of integration.\\
     Substituting back the value of $\sec\theta$ and $\tan\theta$ with respect to $t$ we have,
     \begin{equation}
     \begin{split}
     q_{2}=C_{1}\left[\ln\left(\frac{t+C_{2}+\sqrt{(t+C_{2})^2+C_{1}^2}}{C_{1}}\right)\right]+C_{3}
     \end{split}
     \end{equation}
      Now from equations (\ref{momentum wrt position}) we have,
      \begin{equation}
      p_{1}=\dot{q}_{1}=\frac{t+C_{2}}{\sqrt{(t+C_{2})^2+C_{1}^2}}
      \end{equation}
      and,
       \begin{equation}
      p_{2}=\dot{q}_{2}=\frac{C_{1}}{\sqrt{(t+C_{2})^2+C_{1}^2}}
      \end{equation}
      Therefore we have the position momentum $(q,p)$ set as,
      \begin{equation}
      \begin{split}
      & q_{1}=\sqrt{(t+C_{2})^2+C_{1}^2} \\ &
      q_{2}=C_{1}\left[\ln\left(\frac{t+C_{2}+\sqrt{(t+C_{2})^2+C_{1}^2}}{C_{1}}\right)\right]+C_{3} \\ &
      p_{1}=\frac{t+C_{2}}{\sqrt{(t+C_{2})^2+C_{1}^2}} \\ &
      p_{2}=\frac{C_{1}}{\sqrt{(t+C_{2})^2+C_{1}^2}}
      \end{split}
      \end{equation}
      \section{Computaion of Optimal Controls}
      We have from equations (\ref{Optimal Contols}),
     \begin{equation}
     \begin{split}
     & u_{1}=<p,f_{1}(q)> \\ &
     u_{2}=<p,f_{2}(q)>
     \end{split}
     \end{equation}
     The $1st$ equation of the above set implies
     \begin{equation}
     \begin{split}
     u_{1}=p^{T}f_{1}(q)& =\begin{bmatrix}
     p_{1} & p_{2}
     \end{bmatrix}\begin{bmatrix}
     q_{1} \\ q_{2}
     \end{bmatrix} =p_{1}q_{1}+p_{2}q_{2}\\&=t+C_{2}+\frac{C_{1}}{\sqrt{(t+C_{2})^2+C_{1}^2}}\left[C_{1}\left[\ln\left(\frac{t+C_{2}+\sqrt{(t+C_{2})^2+C_{1}^2}}{C_{1}}\right)\right]+C_{3}\right]
     \end{split}
     \end{equation}
     The $2nd$ equation of the above set implies
     \begin{equation}
     \begin{split}
     u_{2}=p^{T}f_{2}(q)& =\begin{bmatrix}
     p_{1} & p_{2}
     \end{bmatrix}\begin{bmatrix}
     q_{2} \\ -q_{1}
     \end{bmatrix} = p_{1}q_{2}-p_{2}q_{1} \\ &
     =\frac{t+C_{2}}{\sqrt{(t+C_{2})^2+C_{1}^2}}\left[C_{1}\left[\ln\left(\frac{t+C_{2}+\sqrt{(t+C_{2})^2+C_{1}^2}}{C_{1}}\right)\right]+C_{3}\right]-C_{1}
     \end{split}
     \end{equation}
     Now we show the modulus of the control is 1. For that we have,
     \begin{equation}
     \begin{split}
     &u_{1}^2+u_{2}^2 \\ &
     =(t+C_{2})^2+\frac{C_{1}^2}{(t+C_{2})^2+C_{1}^2}\left[C_{1}\left[\ln\left(\frac{t+C_{2}+\sqrt{(t+C_{2})^2+C_{1}^2}}{C_{1}}\right)\right]+C_{3}\right]^2 \\ &+\frac{(t+C_{2})^2}{(t+C_{2})^2+C_{1}^2}\left[C_{1}\left[\ln\left(\frac{t+C_{2}+\sqrt{(t+C_{2})^2+C_{1}^2}}{C_{1}}\right)\right]+C_{3}\right]^2+C_{1}^2 \\ &
     =(t+C_{2})^2+C_{1}^2+\left[C_{1}\left[\ln\left(\frac{t+C_{2}+\sqrt{(t+C_{2})^2+C_{1}^2}}{C_{1}}\right)\right]+C_{3}\right]^2 \\ &
     =q_{1}^2+q_{2}^2 \\ &
     =1 \; \; [\text{by equation \ref{position eqn}}]
     \end{split}
     \end{equation}
     This shows that the norm of the control is $1$ which should always be the case when the normal extremals are arc length parametrized.
     \section{Computation Of the Curvatures}
     We know $u_{1}^2+u_{2}^2=1$. So we can define $u_{1}=\cos(\theta)$ and $u_{2}=\sin(\theta)$ for some $\theta$. Then considering a $(q,\theta)$ coordinate the following equations can be written down as shown in \cite{agrachev2019comprehensive}
     \begin{equation}\label{q theta coordinate}
     \begin{split}
     & \dot{\theta}=c_{1}(q)\cos(\theta)+c_{2}(q)\sin(\theta) \\ &
     \dot{q}=\cos(\theta)f_{1}(q)+\sin(\theta)f_{2}(q)
     \end{split}
     \end{equation}
     where $c_{1},c_{2} \in C^{\infty}(M)$ such that\cite{agrachev2019comprehensive}
     \begin{equation}\label{Lie Bracket}
     [f_{1},f_{2}]=c_{1}f_{1}+c_{2}f_{2}
     \end{equation}
     This suggests that an arc length parametrized normal extremal(that is which satifies the second equation of the set (\ref{q theta coordinate})) satisfies the first equation of the set (\ref{q theta coordinate}). This suggests that the geodesic curvature of the trajectory on $M$ can be written as\cite{agrachev2019comprehensive}
     \begin{equation}\label{formula geodesic}
     \kappa_{g}=\dot{\theta}-c_{1}(q)\cos(\theta)-c_{2}(q)\sin(\theta)
     \end{equation}
     Now we need to compute the functions $c_{1}$ and $c_{2}$. For that we need to compute the Lie bracket $[f_{1},f_{2}]$. If $X_{1},X_{2}$ are the components of $f_{1}$ and $Y_{1},Y_{2}$ are the components of $f_{2}$, denoting $\frac{\partial}{\partial q_{1}}$ by $\partial_{1}$ and $\frac{\partial}{\partial q_{2}}$ by $\partial_{2}$ \\
     first component of the vector $[f_{1},f_{2}]$ is given by
     \begin{equation}
     \begin{split}
     & X^1\partial_{1}Y^1-Y^1\partial_{1}X^1+X^2\partial_{2}Y^1-Y^2\partial_{2}X^1 \\ &
     = 0-q_{2}+q_{2}-0 \\ &
     =0
     \end{split}
     \end{equation}
     second component of the vector $[f_{1},f_{2}]$ is given by
     \begin{equation}
     \begin{split}
     & X^1\partial_{1}Y^2-Y^1\partial_{1}X^2+X^2\partial_{2}Y^2-Y^2\partial_{2}X^2 \\ &
     =-q_{1}-0+0+q_{1} \\ &
     =0
     \end{split}
     \end{equation}
     Therefore from equation (\ref{Lie Bracket}),
     \begin{equation}
     \begin{split}
    & 0=[f_{1},f_{2}]=c_{1}f_{1}+c_{2}f_{2} \\ &
    \implies c_{1}=c_{2}=0
     \end{split}
     \end{equation}
     Now we have considered $u_{1}=\cos(\theta)$. Therefore $\theta=\cos^{-1}(u_{1})$. Therefore,
     \begin{equation}
     \begin{split}
     \theta=\cos^{-1}\left[t+C_{2}+\frac{C_{1}}{\sqrt{(t+C_{2})^2+C_{1}^2}}\left[C_{1}\left[\ln\left(\frac{t+C_{2}+\sqrt{(t+C_{2})^2+C_{1}^2}}{C_{1}}\right)\right]+C_{3}\right]\right]
     \end{split}
     \end{equation} 
     Therefore,
     \begin{equation}
     \begin{split}
     \dot{\theta}=-\frac{1-\frac{C_{1}(t+C_{2})}{[(t+C_{2})^2+C_{1}]^{\frac{3}{2}}}\left[C_{1}\left[\ln\left(\frac{t+C_{2}+\sqrt{(t+C_{2})^2+C_{1}^2}}{C_{1}}\right)\right]+C_{3}\right]+\frac{C_{1}^2}{(t+C_{2})^2+C_{1}^2}}{\sqrt{1-\left[t+C_{2}+\frac{C_{1}}{\sqrt{(t+C_{2})^2+C_{1}^2}}\left[C_{1}\left[\ln\left(\frac{t+C_{2}+\sqrt{(t+C_{2})^2+C_{1}^2}}{C_{1}}\right)\right]+C_{3}\right]\right]^2}}
     \end{split}
     \end{equation}
     Also $c_{1}=c_{2}=0$. \\ Therefore by equation (\ref{formula geodesic})the geodesic curvature, 
     \begin{equation} \kappa_{g}=-\frac{1-\frac{C_{1}(t+C_{2})}{[(t+C_{2})^2+C_{1}]^{\frac{3}{2}}}\left[C_{1}\left[\ln\left(\frac{t+C_{2}+\sqrt{(t+C_{2})^2+C_{1}^2}}{C_{1}}\right)\right]+C_{3}\right]+\frac{C_{1}^2}{(t+C_{2})^2+C_{1}^2}}{\sqrt{1-\left[t+C_{2}+\frac{C_{1}}{\sqrt{(t+C_{2})^2+C_{1}^2}}\left[C_{1}\left[\ln\left(\frac{t+C_{2}+\sqrt{(t+C_{2})^2+C_{1}^2}}{C_{1}}\right)\right]+C_{3}\right]\right]^2}}
     \end{equation}
     Note that the geodesic curvature depends on time $t$. \\
     The Gaussian curvature is given by\cite{agrachev2019comprehensive},
     \begin{equation}
     \kappa=f_{1}(c_{2})-f_{2}(c_{1})-c_{1}^2-c_{2}^2
     \end{equation}
     Now $c_{1}=c_{2}=0$. Therefore $f_{1}(c_{2})=f_{2}(c_{1})=0$ which implies the gaussian curvature,
     \begin{equation}
    \kappa=0
     \end{equation}
     which should be case when $f_{1},f_{2}$ are involutive that is their Lie bracket generates $0$.
     \section{Conclusion}
     We try to conclude the article with a modest conclusion. We want to raise the important points of our work. First thing we want to point out is about the nature of optimality that we have considered in our work. What we have done is actually spacial optimality and not time optimality as time is fixed between $0$ and $1$. In this constraint we obtained the normal extremals exhibited by the trajectory using Pontryagin Minimization formalism. The normal extremals are supposed to be the geodesic on the manifold. The Hamilton's equations that came up in this regard generated nonlinear differential equations of time which are hard to solve. But using the fact that normal extremals are parametrized by arc length which in turn implies the velocity at each point of the trajectory is $1$, we were able to solve the equations successfully. The optimal controls were solved  with respect to time $t$. We further showed that the modulus of the control is $1$ which is always the case in case of arc length parametrized extremals. Further we computed the geodesic curvature and the Gaussian curvature of the Riemannian structure. The geodesic curvature is found to depend on time $t$ while the Gaussian curvature evaluates to $0$. We wish to comment about the result of Gaussian curvature that we found. The reason for that, the orthogonal vector fields which we have considered for our work are found to be involutive that is their Lie Bracket generates $0$. The curvature is always $0$ in the case when the family of vector fields is involutive.
     
     \bibliographystyle{unsrt}
     \bibliography{library}
\end{document}